
\bigskip
\magnification=1200
\parindent=0cm \parskip=3pt plus 1pt minus 1pt
\hsize=13.75cm \vsize=20cm

\font\hugebf=cmbx10 at 14.4pt
\font\bigbf=cmbx10 at 12pt

\font\eightrm=cmr10 at 8pt
\font\eightit=cmti10 at 8pt

\font\tenCal=eusm10
\font\sevenCal=eusm7
\font\fiveCal=eusm5
\newfam\Calfam
  \textfont\Calfam=\tenCal
  \scriptfont\Calfam=\sevenCal
  \scriptscriptfont\Calfam=\fiveCal
\def\mathcal{\fam\Calfam\tenCal}

\font\tenmsa=msam10
\font\sevenmsa=msam7
\font\fivemsa=msam5
\newfam\msafam
  \textfont\msafam=\tenmsa
  \scriptfont\msafam=\sevenmsa
  \scriptscriptfont\msafam=\fivemsa

\font\tenmsb=msbm10
\font\sevenmsb=msbm7
\font\fivemsb=msbm5
\newfam\msbfam
  \textfont\msbfam=\tenmsb
  \scriptfont\msbfam=\sevenmsb
  \scriptscriptfont\msbfam=\fivemsb
\def\mathbb{\fam\msbfam\tenmsb}

\def\bC{{\mathbb C}}
\def\bG{{\mathbb G}}
\def\bK{{\mathbb K}}
\def\bN{{\mathbb N}}
\def\bP{{\mathbb P}}
\def\bQ{{\mathbb Q}}

\def\bZ{{\mathbb Z}}

\def\cF{{\mathcal F}}

\def\cJ{{\mathcal J}}
\def\cK{{\mathcal K}}
\def\cL{{\mathcal L}}
\def\cO{{\mathcal O}}

\def\hexnbr#1{\ifnum#1<10 \number#1\else
 \ifnum#1=10 A\else\ifnum#1=11 B\else\ifnum#1=12 C\else
 \ifnum#1=13 D\else\ifnum#1=14 E\else\ifnum#1=15 F\fi\fi\fi\fi\fi\fi\fi}

\def\msbtype{\hexnbr\msbfam}
\mathchardef\smallsetminus="2\msbtype72   \let\ssm\smallsetminus
\mathchardef\subsetneq="3\msbtype28

\def\GG{{\rm GG}}
\def\FS{{\rm FS}}
\def\Tr{\mathop{\rm Tr}\nolimits}
\def\codim{\mathop{\rm codim}\nolimits}
\def\rank{\mathop{\rm rank}\nolimits}
\def\ord{\mathop{\rm ord}\nolimits}
\def\reg{\mathop{\rm reg}\nolimits}
\def\Bs{\mathop{\rm Bs}\nolimits}
\def\IBs{\mathop{\rm IBs}\nolimits}
\def\Pic{\mathop{\rm Pic}\nolimits}

\let\ul\underline
\def\\{\hfil\break}

\def\today{\ifcase\month\or
January\or February\or March\or April\or May\or June\or July\or August\or
September\or October\or November\or December\fi \space\number\day,
\number\year}

\long\def\claim#1|#2\endclaim{\par\vskip 5pt\noindent 
{\bf #1.}\ {\it #2}\par}

\def\bibitem#1&#2&#3&#4&%
{\hangindent=1.66cm\hangafter=1
\noindent\rlap{\bf #1}\kern1.66cm{\rm #2} {\it #3}, {\rm #4}.}

\def\frac#1#2{{#1\over #2}}
\def\bu{\hbox{$\scriptstyle\bullet$}}

\def\ddbar{\mathop{\partial\overline\partial}}

\def\square{\null\hfill{\hbox{
\vrule height 1.5ex  width 0.1ex  depth 0ex\kern-0.1ex
\vrule height 1.5ex  width 1.5ex  depth -1.4ex\kern-1.5ex
\vrule height 0.1ex  width 1.5ex  depth 0ex\kern-0.1ex
\vrule height 1.5ex  width 0.1ex  depth 0ex}\kern0.5pt}}
\def\qed{\hfill$\square$}

\centerline{\hugebf On the locus of higher order jets of}
\vskip6pt
\centerline{\hugebf entire curves in complex projective varieties}
\bigskip
\centerline{Jean-Pierre Demailly\footnote{${}^*$}{\eightrm
This work is supported by the ERC grant ALKAGE, grant no.\ 670846
from September 2015.}}
\centerline{Institut Fourier, Universit\'e Grenoble Alpes}
\bigskip\bigskip

{\baselineskip=9.6pt\eightrm
{\eightit Abstract}.
For a given complex projective variety, the existence of entire curves
is strongly constrained by the positivity properties of the cotangent
bundle. The Green-Griffiths-Lang conjecture stipulates that entire
curves drawn on a variety of general type should all be contained in a
proper algebraic subvariety. We present here new results on the
existence of differential equations that strongly restrain the locus
of entire curves in the general context of foliated or directed
varieties, under appropriate positivity conditions.\medskip

{\eightit Keywords.} Projective variety, directed variety, entire
curve, jet differential, Green-Griffiths bundle, Semple bundle,
exceptional locus, algebraic differential operator, holomorphic Morse
inequalities.\medskip

{\eightit MSC classification 2020.} 32Q45, 32H30, 14F06
\vskip2cm}

\line{\hfill\it in memory of Professor C.S.\ Seshadri}
\bigskip

{\bigbf 1. Introduction and goals}

Let $X$ be a complex projective manifold, $\dim_\bC X=n$.
Our aim is to study the existence and distribution of entire curves,
namely, of non constant holomorphic curves $f:\bC\to X$. The global
geometry of $X$ plays a fundamental role in this context, and especially
the positivity properties of the canonical bundle $K_X=\Lambda^nT^*_X$.
One of the major open problems of the domain is the following conjecture
due to Green-Griffiths [GGr80] and Lang [Lan87].

\claim 1.1. GGL conjecture|Assume that $X$ is of general type, namely
that $\kappa(X)=\dim X$ where
$\kappa(X):=\limsup_{m\to+\infty}\log h^0(X,K_X^{\otimes m})/\log m$.
Then there exists an algebraic subvariety $Y\subsetneq X$
containing all entire curves $f:\bC\to X$. 
\endclaim

\claim 1.2. Definition|The smallest algebraic subvariety above
will be denoted $Y={\rm Exc}(X)$ and called the
exceptional locus of~$X$.
\endclaim
\medskip

When $X$ is an arithmetic variety, the exceptional locus is expected to
carry a strong arithmetic significance. Especially, one can
(very optimistically) hope for the following result, which is a slight
variation of a conjecture made by Lang [Lan86]: for a projective
variety $X$ defined over a number field $\bK_0$, the exceptional
locus $Y={\rm Exc}(X)$ in the GGL conjecture coincides with the
Mordell locus, where
${\rm Mordell}(X)$ is the smallest complex subvariety $Y$ such that
$X(\bK)\smallsetminus Y$ is finite for all number fields $\bK\supset\bK_0$.

The GGL conjecture unfortunately seems out of reach at this point.
In the present work, we obtain a number of weaker results that still
provide strong restrictions on the distribution of entire curves
in higher order jet bundles. Among these results, we prove for instance the
following statement.

\claim 1.3. Theorem|Let $X$ be a nonsingular projective variety of
general type. Assume that $\Lambda^2T^*_X$ is strongly big on $X$, in
the sense that for a given ample line bundle $A\in\Pic(X)$ the
symmetric powers $S^m(\Lambda^2T^*_X)\otimes\cO_X(-A)$ are generated
by their global sections on a Zariski open set $X\ssm Y$,
$Y\subsetneq X$, when $m>0$ is large $($if $\Lambda^2T^*_X$ is ample,
we can take $Y=\emptyset)$. Then there exist finitely many
rank~$1$ foliations $\cF_\alpha$ on subvarieties $Z_\alpha\subset X_k$
of a suitable $k$-jet Semple bundle of $X$, such that all entire
curves $f:\bC\to X$ are either contained in $Y$ or have a $k$-jet lifting
$f_{[k]}$ that is contained in some $Z_\alpha$ and tangent
to~$\cF_\alpha$. In particular, the latter curves are supported by
the parabolic leaves of these foliations, which can be parametrized as
a subspace of a finite dimensional variety. 
\endclaim
\medskip

Theorem 1.3 generalizes a result that has been known for a long
time for surfaces of general type, which obviously satisfy the hypotheses
(see [GGr80]). By the work of Etesse [Ete19], another class of examples
of projective manifolds possessing ample exterior powers $\Lambda^pT^*_X$
are general complete intersections $X=H_1\cap\ldots\cap H_c$ of sufficient
high degree in complex projective space $\bP^N_\bC$, when the codimension
$c$ is at least equal to $N/(p+1)$ ([Ete19] provides an explicit bound for
the required degree of the $H_j$'s).

The locus $Y$ described in Theorem certainly
includes all abelian and rationally connected subvarieties $Y\subset X$,
and in these cases, the space of entire curves is infinite dimensional
as soon as $\dim Y\geq 2$, even modulo reparametrization.
Our approach is based on existence theorems for jet differentials,
using holomorphic Morse inequalities ([Dem85], [Dem11]), and involves a finer
study of the geometry of jet bundles. Especially, the proof of Theorem~1.3
relies on the use of certain new tautological morphisms related to
induced directed structures on subvarieties of the higher jet bundles.

A complete solution of the GGL conjecture still appears rather elusive.
The techniques used in sections 6 and 7 suggest that one should try to
make a better use of the \hbox{(semi-)}stability properties of the cotangent
bundle, possibly in connection with Ahlfors currents and related
versions of the vanishing theorem, following for instance the
ideas of McQuillan [McQ98]. We would like to celebrate here the pioneering
work of Professor C.S.\ Seshadri in the study of the positivity and
stability properties of vector bundles, which underlie much of our
approach.
\bigskip

{\bigbf 2. Category of directed varieties}

We are interested in entire curves $f:\bC\to X$ such that $f'(\bC)\subset V$,
where $V$ is a (possibly singular) linear subspace of~$X$, i.e.\
a closed irreducible analytic subspace such that the fiber
$V_x:=V\cap T_{X,x}$ is a vector subspace of the tangent space $T_{X,x}$
for all $x\in X$. We briefly recall
below some of the relevant concepts, and refer to [Dem11], [Dem20] for
further details.

\claim 2.1. Definition of the category of directed varieties|
{\parindent=6.4mm
\item{\rm(a)} Objects are pairs $(X,V)$ where $X$ is a complex manifold and 
$V\subset T_X$ a linear subspace of~$\,T_X$.
\item{\rm(b)} Arrows $\psi:(X,V)\to(Y,W)$ are holomorphic maps $X\to Y$
such that $d\psi(V)\subset W$

\item{\rm(c)} The {\rm absolute case} refers to the case $V=T_X$, i.e.\ of
pairs $(X,T_X)$.
\item{\rm(d)} The {\rm relative case} refers to pairs $(X,T_{X/S})$ where
$X\to S$ is a fibration.
\item{\rm(e)} We say that we are in the {\rm Integrable case} the sheaf of
sections $\cO(V)$ is stable by Lie brackets. This corresponds to 
holomorphic $($and possibly singular$)$ foliations.\vskip0pt}
\endclaim
\medskip

We now define the canonical sheaf of a directed manifold $(X,V)$. When $V$ is
nonsingular, i.e.\ is a subbundle, we simply set $K_V=\det(V^*)$ (this is
a line bundle, i.e.\ an invertible sheaf).
When $V$ is singular, we first introduce the rank $1$ sheaf ${}^b\cK_V$
of sections of $\det V^*$ that are  locally bounded with respect 
to a smooth ambient metric on~$T_X$.\ One can show that ${}^b\cK_V$ is
equal to the integral closure of the image of the natural morphism
$$
\cO(\Lambda^rT_X^*)\to \cO(\Lambda^r V^*)\to \cL_V:=
\hbox{invertible~sheaf}~\cO(\Lambda^r V^*)^{**}
$$
that is, if the image is $\cL_V\otimes\cJ_V$, $\cJ_V\subset\cO_X$,
$$
{}^b\cK_V=\cL_V\otimes\overline{\cJ}_V,~~~~
\overline{\cJ}_V=\hbox{integral closure of}~\cJ_V.
\leqno(2.2)
$$
However, one may have to first blow up $X$ as follows to ensure
that ${}^b\cK_V$ provides the appropriate geometric information.

\claim 2.3. Blow up process for a directed variety|
If $\mu:\widetilde X\to X$ is a modification, then
$\widetilde X$ is equipped with the pull-back directed structure
$\widetilde V=\overline{\tilde\mu^{-1}(V_{|X'})}$, where $X'\subset X$
is a Zariski open set over which $\mu$ is a biholomorphism.
\endclaim

\claim 2.4. Observation|One always has
${}^b\cK_V\subset\mu_*({}^b\cK_{\widetilde V})\subset \cL_V=
\cO(\det V^*)^{**}$,
and $\mu_*({}^b\cK_{\widetilde V})$ ``increases'' with $\mu$~ $($taking
successive blow-ups
$\kern4pt\widetilde{\strut}\kern-5.5pt\widetilde X\to \widetilde X\to X)$.
\endclaim

By Noetherianity, one can define a sequence of rank $1$ sheaves
$$
\cK^{[m]}_V=\lim_{\mu}\uparrow \mu_*({}^b\cK_{\widetilde V})^{\otimes m},~~~
\mu_*({}^b\cK_V)^{\otimes m}\subset\cK^{[m]}_V\subset \cL_V^{\otimes m}
\leqno(2.5)
$$
which we call the {\it pluricanonical sheaf sequence} of $(X,V)$.
Remark that the blow-up $\mu$ for which the limit is attained
may depend on~$m$. We do not know if there exists a modification $\mu$
that works for all~$m$. This generalizes the concept
of {\it reduced singularities} of foliations
(which is known to work in that form only for surfaces).

\claim 2.6. Definition|We say that $(X,V)$ is of {\rm general type} if
the pluricanonical sheaf sequence $\cK_V^{[\bullet]}$ is big, that is, if
$H^0(X,\cK^{[m]}_V)$ provides a generic embedding of $X$ for sufficiently
large powers $m\gg 1$.
\endclaim
\medskip

The Green-Griffiths-Lang conjecture can be generalized to directed varieties
as follows.

\claim 2.7. Generalized GGL conjecture| If $(X,V)$ is a directed manifold of
general type, i.e.\ if $\smash{\cK_V^{[\bullet]}}$ is big,
there exists an algebraic locus $Y\subsetneq X$ such that for every
entire curve $f:(\bC,T_\bC)\to(X,V)$, one has $f(\bC)\subset Y$.
\endclaim
\medskip

When $r=\rank V=1$, the generalized GGL conjecture is an elementary
consequence of the Ahlfors-Schwarz lemma. In fact, the function
$t\mapsto \log\Vert f'(t)\Vert_{V,h}$ is strictly subharmonic
whenever $(V^*,h^*)$ is assumed to be big.\medskip

{\bf 2.8. Remark.} The directed form of the GGL conjecture as stated
above is possibly too optimistic. It~might be safer to add a suitable
$($semi-$)$stability condition on~$V$. In the absolute case $V=T_X$,
such a semistability property is automatically satisfied when $K_X$
is ample, as a consequence of the existence of a K\"ahler-Einstein metric.
\bigskip

{\bigbf 3. Definition of algebraic differential operators.}

The basic strategy to attack the Green-Griffiths-Lang conjecture is to
show that all entire curves $f:(\bC,T_\bC)\to(X,V)$ must satisfy nontrivial
algebraic differential equations $P(f\,;\,f',f'',\ldots,f^{(k)})=0$,
and actually, as we will see, many such
equations. Following [GGr80] and [Dem97], we now introduce the
useful concept of jet differential operator. 
Let $(\bC,0)\to (X,V),~~t\mapsto f(t)=(f_1(t),\ldots,f_n(t))$
be a germ of curve such that $f(0)=x\in X$. Pick local holomorphic
coordinates $(z_1,\ldots,z_n)$ centered at $x$ on a coordinate open set
$U\simeq U'\times U''\subset\bC^r\times\bC^{n-r}$
such that $\pi':U\to U'$ induces
an isomorphism $d\pi':V\to U\times\bC^r$.\ Then $f$ is determined
by the Taylor expansion
$$
\pi'\circ f(t)=t\xi_1+\ldots+t^k\xi_k+O(t^{k+1}),~~~
\xi_s=\frac{1}{s!}\nabla^sf(0),
\leqno(3.1)
$$
where $\nabla$ is the trivial connection on $V\simeq U\times\bC^r$.
One considers the Green-Griffiths bundle $E^\GG_{k,m}V^*$ of polynomials
of weighted degree $m$, written locally in coordinate charts as
$$
P(x\,;\,\xi_1,\ldots,\xi_k)=\sum
a_{\alpha_1\alpha_2\ldots\alpha_k}(x)\xi_1^{\alpha_1}\ldots\xi_k^{\alpha_k},~~~
\xi_s\in V.
$$
These can also be viewed as algebraic differential operators
$$
P(f_{[k]})=P(f\,;\,f',f'',\ldots,f^{(k)})=\sum
a_{\alpha_1\alpha_2\ldots\alpha_k}(f(t))~f'(t)^{\alpha_1}
f''(t)^{\alpha_2}\ldots f^{(k)}(t)^{\alpha_k}.\kern-6pt
\leqno(3.2)
$$
Here $t\mapsto z=f(t)$ is a curve, $f_{[k]}=
(f',f'',\ldots,f^{(k)})$ its $k$-jet, and
$a_{\alpha_1\alpha_2\ldots\alpha_k}(z)$ are supposed to holomorphic 
functions on~$X$. The reparametrization action: $f\mapsto
f\circ\varphi_\lambda$, $\varphi_\lambda(t)=\lambda t$, $\lambda\in\bC^*$ yields
$(f \circ \varphi_\lambda)^{(k)}(t)=\lambda^kf^{(k)}(\lambda t)$,
whence a $\bC^*$-action
$$
\lambda\cdot(\xi_1,\xi_1,\ldots,\xi_k)=(\lambda\xi_1,
\lambda^2\xi_2,\ldots,\lambda^k\xi_k).
\leqno(3.3)
$$
$E^\GG_{k,m}$ is precisely the set of polynomials of weighted
degree $m$, corresponding to coefficients $a_{\alpha_1\ldots\alpha_k}$ with
$m=|\alpha_1|+2|\alpha_2|+\ldots+k|\alpha_k|$. An important fact is the

\claim 3.4. Direct image formula|
Let $J_k^{\rm nc}V\subset J_kV$ be the set of non constant $k$-jets.
One defines the {\rm Green-Griffiths} bundle to be
{\rm $X_k^\GG=J_k^{\rm nc}V/\bC^*$}
and {\rm $\cO_{X_k^\GG}(1)$} to be the associated tautological
rank $1$ sheaf. Let $\pi_k:X_k^\GG\to X$ be the natural projection.
Then we have the direct image formula
$$
E_{k,m}^\GG V^*=(\pi_k)_*\cO_{X_k^\GG}(m).
$$
\endclaim
\medskip

{\bigbf 4. Consequences of the holomorphic Morse inequalities}

Given a real $(1,1)$-form $\theta$ on $X$, the $q$-index set of $\theta$
is defined to be
$$
X(\theta,q)=\{x\in X\,|\;\hbox{$\theta(x)$ has signature $(n-q,q)$}\}
\leqno(4.1)
$$
(exactly $q$ negative eigenvalues and $n-q$ positive ones). We also set
$$
X(\theta,\leq q)=\bigcup_{0\leq j\leq q}X(\theta,j).
\leqno(4.2)
$$
Then $X(\theta,q)$ and $X(\theta,\leq q)$ are open sets.
and sign$(\theta^n)=(-1)^q$ on $X(\theta,q)$. The general statement of
holomorphic Morse inequalities [Dem85] asserts that for any
hermitian line bundle $(L,h)$ of curvature form $\theta=\Theta_{L,h}$
on~$X$, and any coherent sheaf $\cF$ of rank $r=\rank\cF$, one has
$$
\sum_{j=0}^q(-1)^{q-j}h^j(X,L^{\otimes m}\otimes\cF)\leq
r\,{m^n\over n!}\int_{X(\theta,\leq q)}(-1)^q\theta^n+o(m^n)
\quad\hbox{ as $m\to+\infty$}.
\leqno(4.3)
$$
An application of these inequalities yields the following fundamental
estimates~[Dem11].

\claim 4.4. Main cohomology estimates|
Let $(X,V)$ be a directed manifold, $A$ an ample
line bundle over $X$, $(V,h)$ and $(A,h_A)$ hermitian
structures such that $\Theta_{A,h_A}>0$. Define
$$
\leqalignno{
&L_{k,\varepsilon}=\cO_{X_k^\GG}(1)\otimes\pi_k^*\cO\Big
(-{1\over kr}\Big(1+{1\over 2}+\ldots+{1\over k}\Big)\varepsilon A\Big),
&\hbox{\rm(a)}\cr
&\eta_\varepsilon=\Theta_{\det V^*,\det h^*}-\varepsilon\,\Theta_{A,h_A}.
&\hbox{\rm(b)}\cr}
$$
Then all $m\gg k\gg 1$ such that $m$ is sufficiently divisible
and for all $q\ge 0$ we have upper and lower bounds
$$
\leqalignno{
&h^q(X_k^\GG,\cO(L_{k,\varepsilon}^{\otimes m}))\le{m^{n+kr-1}\over (n{+}kr{-}1)!}
{(\log k)^n\over n!\,(k!)^r}\bigg(
\int_{X(\eta,q)}\!\!\!(-1)^q\eta_\varepsilon^n+\frac{C}{\log k}\bigg),
&\hbox{\rm(c)}\cr
&h^q(X_k^\GG,\cO(L_{k,\varepsilon}^{\otimes m}))\ge{m^{n+kr-1}\over (n{+}kr{-}1)!}
{(\log k)^n\over n!\,(k!)^r}\bigg(
\int_{X(\eta,\{q,\,q\pm 1\})}(-1)^q\eta_\varepsilon^n-\frac{C}{\log k}\bigg).
&\hbox{\rm(d)}\cr}
$$
\endclaim

The case $q=0$ is the most useful one, as it gives estimates for the
number of holomorphic sections. 
We now explain the essential ideas involved in the proof of the estimates.
\medskip

{\it $1{}^{\rm st}$ step: construction of a Finsler metric on $k$-jet bundles}

Let $J_kV$ be the bundle of $k$-jets of curves $f:(\bC,T_\bC)\to(X,V)$.
Assuming that $V$ is equipped with a hermitian metric $h$,
one defines a ''weighted Finsler metric'' on $J^kV$ by taking $p=k!$ and
$$
\Psi_{h_k}(f):=\Big(\sum_{1\le s\le k}\big(\varepsilon_s\Vert\nabla^sf(0)
\Vert_{h(x)}\big)^{2p/s}\Big)^{1/p},~~1=\varepsilon_1\gg
\varepsilon_2\gg\cdots\gg\varepsilon_k.
\leqno(4.5)
$$
Letting $\xi_s=\nabla^sf(0)$, this can actually be viewed as a 
metric $h_k$ on $L_k:=\cO_{X_k^\GG}(1)$. Some error terms of
the form $O(\varepsilon_s/\varepsilon_t)_{t<s}$ appear, but they are
negligible in the limit, and the leading term of the curvature form is
$(z,\xi_1,\ldots,\xi_k)\mapsto \Theta^0_{L_k,h_k}(z,\xi)$ with
$$
\Theta^0_{L_k,h_k}(z,\xi)=\omega_{\FS,k}(\xi)+{i\over 2\pi}
\sum_{1\le s\le k}{1\over s}{|\xi_s|^{2p/s}\over \sum_t|\xi_t|^{2p/t}}
\sum_{i,j,\alpha,\beta}c_{ij\alpha\beta}(z)
{\xi_{s\alpha}\overline\xi_{s\beta}\over|\xi_s|^2}\,dz_i\wedge d\overline z_j,
\leqno(4.6)
$$
where $(c_{ij\alpha\beta}(s))$ are the coefficients of the curvature tensor
$\Theta_{V^*,h^*}$ at point $z$, and $\omega_{\FS,k}(\xi)$ is the vertical
Fubini-Study metric on the fibers of $X_k^\GG\to X$.
The expression gets simpler by using polar coordinates
$$
x_s={\vert\xi_s\vert_h^{2p/s}\over \sum_t\vert\xi_t\vert_h^{2p/t}},\quad
u_s={\xi_s\over \vert\xi_s\vert_h}={\nabla^sf(0)\over\vert\nabla^sf(0)\vert},
\quad 1\leq s\leq k.
\leqno(4.7)
$$

{\it $2{}^{\rm nd}$ step: probabilistic interpretation of the curvature}

In such polar coordinates, one gets the formula
$$
\Theta_{L_k,h_k}(z,\xi)\simeq\omega_{\FS,p,k}(\xi)+
{i\over 2\pi}\sum_{1\le s\le k}{1\over s}x_s
\sum_{i,j,\alpha,\beta}c_{ij\alpha\beta}(z)\,
u_{s\alpha}\overline u_{s\beta}\,dz_i\wedge d\overline z_j
\leqno(4.8)
$$
where $\omega_{\FS,k}(\xi)$ is positive definite in $\xi$. The other 
terms are a weighted average of the values of the 
curvature tensor $\Theta_{V,h}$ on vectors $u_s$ in the unit sphere
bundle $SV\subset V$. The weighted projective space can be viewed
as a circle quotient of the pseudosphere $\sum|\xi_s|^{2p/s}=1$, so
we can take here $x_s\ge 0$, $\sum x_s=1$. This is essentially a
sum of the form $\sum\frac{1}{s}Q(u_s)$ where
$Q(u)=\langle\Theta_{V^*,h^*}u,u\rangle$ and 
$u_s$ are random points of the sphere, and so as $k\to+\infty$ this
can be estimated by a ``Monte-Carlo'' integral
$$
\Big(1+\frac{1}{2}+\ldots+\frac{1}{k}\Big)\int_{u\in SV}Q(u)\,du.
\leqno(4.9)
$$
As $Q$ is quadratic, we have
$$\int_{u\in SV}
Q(u)\,du=\frac{1}{r}\Tr(Q)\,{=}\,
\frac{1}{r}\Tr(\Theta_{V^*,h^*})\,{=}\,
{1\over r}\Theta_{\det V^*,\det h^*}.
\leqno(4.10)
$$
The above equality show that the Monte Carlo approximation only depends
on $\det V^*$ and that for the unitary invariant probability measure
taken on $(SV)^k\ni(u_s)$, we have an expected value
$$
{\bf E}(\Theta_{L_k,h_k})=
\omega_{\FS,p,k}(\xi)+\sum_{1\leq s\leq k}{1\over s}
{1\over kr}\Theta_{\det V^*,\det h^*}(z)
\leqno(4.11)
$$
(The factor ${1\over k}$ comes from the fact that the expected value
of $x_s$ on the $(k-1)$-simplex $\sum_{1\leq s\leq k}x_s=1$ is ${1\over k}$).
Formula (4.11) is the main reason why the leading term of the
cohomology estimates only involves $\det V^*$. Of course, a complete
proof requires an estimate of the standard deviation occurring in the
Monte Carlo process, as is done in [Dem11].\qed
\medskip

By passing to a ``singular version'' of holomorphic Morse inequalities
to accommodate singular metrics (Bonavero, [Bon93]) and subtracting
the terms involving $\varepsilon\,\Theta_{A,h_A}$, one~gets

\claim 4.12. Corollary: existence of global jet differentials {\rm ([Dem11)}|
Assume that $(X,V)$ is of general type, namely that ${}^b\cK_V$
is a {\rm big rank $1$ sheaf}, and let
$$
L_{k,\varepsilon}=\cO_{X_k^\GG}(1)\otimes
\pi_k^*\cO(-\delta_k\varepsilon A), \quad
\delta_k={1\over kr}\Big(1+{1\over 2}+\cdots+{1\over k}\Big)
$$
with $A\in\Pic(X)$ ample. Then there exist many nontrivial global sections
$$
P\in H^0(X_k^\GG,L_{k,\varepsilon}^{\otimes m})\simeq
H^0(X,E^\GG_{k,m}V^*\otimes\cO(-m\delta_k\varepsilon A))    
$$
for $m\gg k\gg 1$ and $\varepsilon\in\bQ_{>0}$ small.
\endclaim
\medskip

The fact that entire curves satisfy differential equations is now a
consequence of the following result.

\vbox{\claim 4.13. Fundamental vanishing theorem {\rm ([GGr80],
[Dem97], [SYe97])}|For all global differential operators
$$
P\in H^0(X,E^\GG_{k,m}V^*\otimes\cO(-qA)),\quad q\in\bN^*,
$$
and all $f:(\bC,T_\bC)\to (X,V)$, one has $P(f_{[k]})\equiv 0$.
\endclaim}
\medskip

Geometrically, this can interpreted by stating that the image $f_{[k]}(\bC)$ of the $k$-jet curve lies in the base
locus
$$
Z=\bigcap_{m\in\bN^*}~
\bigcap_{\sigma\in H^0(X_k^\GG,L_{k,\varepsilon}^{\otimes m})}
\sigma^{-1}(0)~~\subset~~X^\GG_k.
\leqno(4.14)
$$
To prove the GGL conjecture, we would need to check that $\pi_k(Z)\subsetneq X$.
This turns out to be a hard problem.
\bigskip

{\bigbf 5. Investigation of the base locus}

We start by formulating the base locus problem in a very general context.
Let $(L,h)$ be a hermitian line bundle over $X$. If we assume
that $\theta=\Theta_{L,h}$ satisfies $\int_{X(\theta,\leq 1)}\theta^n>0$,
then we know that $L$ is big, i.e.\ that
$h^0(X,L^{\otimes m})\geq c\,m^n$, for $m\geq m_0$ and $c>0$,
but this does not tell us anything about the base locus
$\Bs(L)=\bigcap_{\sigma\in H^0(X,L^{\otimes m})}\sigma^{-1}(0)$.

\claim 5.1. Definition|
The ``iterated base locus'' $\IBs(L)$ is obtained by picking inductively
$Z_0=X$ and $Z_k={}$zero divisor of a section $\sigma_k$ of $L^{\otimes m_k}$
over the normalization of $Z_{k-1}$, and taking
$\bigcap_{k,m_1,\ldots,m_k,\sigma_1,\ldots,\sigma_k}Z_k$.
\endclaim

\claim 5.2. Unsolved problem|
Find a condition, e.g.\ in the form of Morse integrals $($or analogous
integrals$)$ for $\theta=\Theta_{L,h}$, ensuring that $\codim\IBs(L)>p$.
\endclaim
\medskip

For instance, there might exist a way of deriving from such conditions
the positivity of Morse integrals $\int_{Z(\theta_{|Z},\leq 1)}\theta^{n-p}$
for arbitrary irreducible subvarieties, $\codim Z=p$, that are obtained
themselves as iterated base loci. In the specific case of the
Green-Griffiths tautological line bundle $L_k$, one can get the following
statement (without loss of generality, we can exclude the case of directed
structures of rank $r=1$, since the GGL conjecture is then trivial).

\claim 5.3. Theorem|Let $(X,V)$ be a directed variety of {\rm general type},
with $r=\rank V\geq 2$. Then there exist $k_0\in\bN$
and $\delta>0$ with the following properties.
Let $Z\subset X^\GG_k$ be an algebraic subvariety that is
a complete intersection of irreducible hypersurfaces
$$
Z=\bigcap_{1\le j\le \ell}\big\{\hbox{$k$-jets}~f_{[k]}\in X^\GG_k\,;\,P_j(f)=0
\big\},~~P_j\in H^0(X,E^\GG_{s_j,m_j}V^*\otimes G_j),
$$
with $k\geq k_0$, $\ord(P_j)=s_j$, $1\leq s_1<\cdots<s_\ell\leq k$,
$\sum\limits_{1\leq j\leq\ell} {1\over s_j}\leq \delta\log k$,
and $G_j\in\Pic(X)$. Then $\dim Z=n+kr-1-\ell$ and the Morse integrals
$$
\int_{Z(L_{k,\varepsilon},\leq 1)}\Theta_{L_{k,\varepsilon}}^{\dim Z}
$$
of $L_{k,\varepsilon}=\cO_{X_k^\GG}(1)\otimes
\pi_k^*\cO_X\big(-{1\over kr}\big(1+{1\over 2}+\cdots+{1\over k}\big)
\varepsilon A\big)$
are positive for $\varepsilon>0$ small, hence
$H^0(Z,L_{k,\varepsilon}^{\otimes m})\ge c\,m^{\dim Z}$ for $m\gg 1$, with $c>0$.
\endclaim
\medskip

At the expense of a slightly more involved statement, it would be possible
to allow some repetitions in the degrees $s_j$ of the polynomials $P_j$.
Unfortunately, none of the extended versions we have been able to reach
seems sufficient to prove the GGL conjecture, which would follow if we could
cut down $Z$ to a subvariety of dimension $\dim Z<n$, thus of very high
codimension in the jet bundle.
\medskip

{\it Proof}. We come back to the calculations made in [Dem11],
especially those detailed in~\S$\,$2. The main point is an integration
on the fiber of $X_k^\GG\to X$, which is reduced in polar coordinates
to an integration on a product of unit spheres $(SV)^k$ with the
$(k-1)$-dimensional simplex. Here, we
integrate on the fibers $Z_z$ of $Z\to X$, namely on a rather complicated
subvariety of the weighted projective space
$P_{\bu}(V^k)=(V^k\ssm\{0\})/\bC^*$ of the form
$$
P_j(z,\xi_1,\ldots,\xi_{s_j})=0,\quad\hbox{($z\in X$ fixed)},
\leqno(5.4)
$$
where the $P_j$'s are weighted homogeneous polynomials in
$(\xi_1,\ldots,\xi_{s_j})$ of positive degree in $\xi_{s_j}$.
The choice of our Finsler metric can be combined with a rescaling of the form
$\xi_s=\varepsilon_s^{-1}\tilde\xi_s$ (see [Dem11, Lemma 2.12]). Since
$\varepsilon_s\ll\varepsilon_t$ for $t<s$, the rescaled subvariety (5.4)
actually converges to the subvariety defined
in the new coordinates $\tilde\xi=(\tilde\xi_s)$ by the equations
$$
Q_j(z,\tilde\xi_{s_j})=
P_j(z,0,\ldots,0,\tilde\xi_{s_j})=0,
\leqno(5.5)
$$
at least, for all $x$ in the full measure Zariski open set of $X$
where the leading coefficients of $Q_j(z,\bu)$ do not vanish.
These polynomials depend on non overlapping variables $\xi_{s_j}$,
hence the limit subvariety can be seen as a product of cones
over a product subvariety in $P(V)^k$, where some of the factors $P(V)$ are
replaced by hypersurfaces. The calculation of the Morse integrals on $Z$
requires computing the integrals
$$
\int_{Z_z}x^\alpha\omega_{\FS,p,k}(\xi)^N,\quad
x_s={\vert\xi_s\vert_h^{2p/s}\over \sum_t\vert\xi_t\vert_h^{2p/t}},
\quad N=\dim Z_z=kr-1-\ell,\quad|\alpha|\leq n,
$$
and an evaluation of the standard deviation via the Cauchy-Schwarz inequality
involves the same type of integrals with $|\alpha|\leq 2n-2$ (see [Dem11]).
By the Lelong-Poincar\'e formula and the Fubini theorem, the above integrals
are equal to
$$
\int_{P_{\bu}(V^k)}x^\alpha\omega_{\FS,p,k}(\xi)^N\wedge
\bigwedge_{1\leq j\leq\ell}{i\over 2\pi}\ddbar\log|Q_j(z,\xi_{s_j})|^2.
$$
In view of the unitary invariance of $x^\alpha\omega_{\FS,p,k}(\xi)^N$
in each $\xi_s$, the Crofton formula implies that the previous integral is
equal to
$$
\int_{P_{\bu}(V^k)}x^\alpha\omega_{\FS,p,k}(\xi)^N\wedge
\bigwedge_{1\leq j\leq\ell}d_j\,{i\over 2\pi}\ddbar\log|\xi_{s_j}|_h^2
\leqno(5.6)
$$
where $d_j=\deg_{\xi_{s_j}}Q_j(z,\xi_{s_j})$. Here the Fubini-Study metric
is derived from the Finsler metric $(\sum_s\vert\xi_s\vert_h^{2p/s})^{1/p}$.
If we make a change of variable $\xi_{s,\lambda}\mapsto \xi_{s,\lambda}^s$
for each component of $\xi_s$ and apply the comparison principle of [Dem87],
we see that (5.6) is equivalent, up to factors controlled by constants
$C_n^{\pm 1}$ depending only on the dimension of $X$, to the integral
$$
{1\over(k!)^r}\prod_{1\leq j\leq\ell}d_js_j
\int_{P(\bC^{kr})}x^\alpha\omega_\FS(\xi)^N\wedge
\bigwedge_{1\leq j\leq\ell}{i\over 2\pi}\ddbar\log|\xi_{s_j}|^2
\leqno(5.7)
$$
where $|\xi|^2=\sum_s|\xi_s|^2$,
$\omega_\FS(\xi)={i\over 2\pi}\ddbar\log|\xi|^2$ and
$x_s={|\xi_s|^2 \over|\xi|^2}$ . Additionally, the correcting factor
is equal to $1$ if $\alpha=0$, and in that case the precise value of
our integrals is ${1\over(k!)^r}\prod_{1\leq j\leq\ell}d_js_j$,
where the factor ${1\over (k!)^r}$ comes from the weighted Fubini-Study
metric, and $\prod_{1\leq j\leq\ell}d_js_j$ is the relative degree
of $Z_z$ in the weighted projective space. Again by the Crofton formula,
we can replace each factor ${i\over 2\pi}\ddbar\log|\xi_{s_j}|^2$
by the current of integration on the hyperplane  $\xi_{s_j,1}=0$. This shows
that our integrals are equal to
$$
{1\over(k!)^r}\prod_{1\leq j\leq\ell}d_js_j
\int_{P((\bC^r)^{k-\ell}\times(\bC^{r-1})^\ell)}x^\alpha\omega_\FS(\xi)^N,
\leqno(5.8)
$$
and these can be computed by means of the formulas given in [Dem11].
In particular, the calculation
of the expected value ${\bf E}(\Theta_{L_k,h_k})$ depends on the
integrals (5.7) with $|\alpha|=1$, i.e.\ $x^\alpha=x_s$. There are
2 possible values, one for $s\in\{s_1,\ldots,s_\ell\}$, namely
$I'={r-1\over kr-\ell}$ and another for $s$ taken in the complement,
equal to $I''={r\over kr-\ell}$, so that $\ell I'+(k-\ell)I''=1$.
The ratio $I'/I''$ belongs to $[1-1/r,1]\subset[1/2,1]$.
Our assumption $\sum_{1\leq j\leq\ell}{1\over s_j}\leq\delta\log k$
implies $\log\ell=\delta\log k+O(1)$, hence
$\ell=O(k^\delta)\ll k$. Therefore, if
we are concerned with estimates only up to universal constants,
we may consider $I'$, $I''$ to be of the same order of magnitude ${1\over k}$.
Then, most of the arguments employed in [Dem11] remain unchanged. The
integration is performed on a certain codimension $\ell$ subvariety
of $\bC^{kr}$ which is a ramified cover over a product
$\bC^{kr-\ell}=(\bC^r)^{k-\ell}\times(\bC^{r-1})^\ell$, and the equality
(4.10) applies to each of the $(k-\ell)$ factors $\bC^r$, where the polar
coordinate change yields an integral on the whole unit sphere. The $\ell$
remaining factors possibly lead to terms of unknown signature, but their sum
has a relative size $O(\sum{1\over s_j}/\sum{1\over s})=O(\delta)$.
Therefore the expected value of the curvature
${\bf E}(\Theta_{L_{k,\varepsilon},h_k})$, which can be derived from
(4.11) (see also [Dem11, 2.18]), has an horizontal term equal to
${1\over kr}\log k\,(\eta_\varepsilon\pm O(\delta))$.
Our estimate follows.\qed\medskip

{\bf 5.9. Remark.} The above calculations even give an evaluation of the
leading constant $c$, namely
$$
c={\prod_{1\leq j\leq\ell}d_js_j
\over(n+kr-1)!}{(\log k)^n\over n!(k!)^r}\Big(\int_{X(\eta_\varepsilon,\leq 1)}
\eta_\varepsilon^n-O((\log k)^{-1})-O(\delta)\Big)
$$
where $\eta_\varepsilon=\Theta_{\det V^*,\det h^*}-\varepsilon\,\Theta_{A,h_A}$.
If $Z$ is taken to be just one irreducible component of a complete
intersection, the quantity $\prod_{1\leq j\leq\ell}d_js_j$ has to be replaced
by the relative degree over $X$ of that component.
\bigskip

{\bigbf 6. Semple jet bundles and induced directed structures}

Semple jet bundles provide further geometric information that appear very
useful in the investigation of the base locus. We briefly recall the relevant
notation and concepts from [Dem97].  To begin with,
we assume that $X$ is non singular and connected, and that $V$ is non
singular, i.e.\ that $V$ is a subbundle of~$T_X$, and we set $\dim X=n$,
$\rank V=r$.

\claim 6.1. 1-jet functor|This is the functor $(X,V)\mapsto (\tilde X,\tilde V)$
on the category of directed varieties defined by
{\parindent=6.5mm
\item{\rm(a)}
$\tilde X=P(V)={}$bundle of projective spaces of lines in $V,$
\item{\rm(b)}
$\pi:\tilde X=P(V)\to X,~~~(x,[v])\mapsto x,~~v\in V_x\,,$
\item{\rm(c)}
$\tilde V_{(x,[v])}=\big\{\xi\in T_{\tilde X,(x,[v])}\,;\;\pi_*\xi\in\bC v \subset T_{X,x}\big\}.$}
\endclaim
\medskip

By taking tangents in $\tilde V=P(V)$, every entire
curve $f:(\bC,T_\bC)\to(X,V)$
lifts as a curve $f_{[1]}:(\bC,T_\bC)\to(\tilde X,\tilde V)$ by putting
$f_{[1]}(t):=(f(t),[f'(t)])\in P(V_{f(t)})\subset \tilde X$.

\claim 6.2. Definition of Semple jet bundles|
{\parindent=6.5mm
\item{\rm(a)} We define $(X_k,V_k)$ to be the $k$-th iteration of the functor 
$(X,V)\mapsto(\tilde X,\tilde V)$.
\item{\rm(b)} In this way, every curve $f:(\bC,T_\bC)\to(X,V)$ gives rise to
a projectivized $k$-jet lifting $f_{[k]}:(\bC,T_\bC)\to(X_k,V_k)$.\vskip0pt}
\endclaim
\medskip

{\bf 6.3. Basic exact sequences.} By construction we have $X_k=P(V_{k-1})$
and a tautological line bundle $\cO_{X_k}(1)$ on $X_k=P(V_{k-1})$.
One can also check that there are exact sequences
$$
\leqalignno{
&0\to T_{X_k/X_{k-1}}\to V_k
\mathop{\longrightarrow}\limits^{d\pi_k}_{}\cO_{X_k}(-1)
\to 0,&\hbox{\rm(a)}\cr
\noalign{\vskip4pt}      
&0\to\cO_{X_k}\to \pi_k^\star V_{k-1}\otimes\cO_{X_k}(1)
\to T_{X_k/X_{k-1}}\to 0.&\hbox{\rm(b)}\cr}
$$
The sequence (a) is equivalent to the definition of $V_k$, and (b) is just
the Euler exact sequence for $X_k=P(V_{k-1})$.
These sequences imply that $\rank V_k=\rank V$ is constant.
Letting $n=\dim X$ and $r=\rank V$, we obtain morphisms
$$
\pi_{k,\ell}=\pi_{\ell-1}\circ\cdots\circ\pi_{k-1}\circ\pi_k:X_k\to X_\ell
$$
and a {\it tower of $\bP^{r-1}$-bundles}
$$
\pi_{k,0}:X_k\mathop{\longrightarrow}\limits^{\pi_k}_{}X_{k-1}
\to\cdots\to X_1\mathop{\longrightarrow}\limits^{\pi_1}_{}X_0=X
\leqno(6.4)
$$
with $\dim X_k=n+k(r-1)$, $\rank V_k=r$. In the sequel, we introduce
the weighted invertible sheaves
$$
\cO_{X_k}({\ul a})=\pi_{k,1}^*\cO_{X_1}(a_1)\otimes{\cdot}{\cdot}{\cdot}
\otimes\pi_{k,k-1}^*\cO_{X_{k-1}}(a_{k-1})\otimes\cO_{X_k}(a_k)
\leqno(6.5)
$$
for~every $k$-tuple ${\ul a}=(a_1,...,a_k)\in\bZ^k$, and let
${\ul 1}=(1,...,1)\in\bZ^k$. Then
$$
\det V_k=\det T_{X_k/X_{k-1}}\otimes\cO_{X_k}(-1),\quad
\det T_{X_k/X_{k-1}}=\pi_k^*\det V_{k-1}\otimes\cO_{X_k}(r),
\leqno(6.6)
$$
and we infer inductively
$$
\leqalignno{
\det V_k&=\pi_k^*\det V_{k-1}\otimes\cO_{X_k}(r-1),&(6.6')\cr
\noalign{\vskip5pt}
&=\pi_{k,0}^\star\det V\otimes\cO_{X_k}((r-1)\ul 1).&(6.6'')\cr}
$$

\claim 6.7. Proposition {\rm ([Dem97])}|When $(X,V)$ is non singular,
the Semple bundle $X_k$ is a smooth compactification of the
quotient $X_k^{\GG,\reg}/\bG_k=J_k^{\GG,\reg}/\bG_k$,
where $\bG_k$ is the group of $k$-jets of germs of biholomorphisms
of $(\bC,0)$, acting on the right by reparametrization:
$(f,\varphi)\mapsto f\circ\varphi$, and where $J_k^{\reg}$ is
the space of $k$-jets of regular curves.
\endclaim
\medskip

In the absolute case $V=T_X$, we obtain what we call the
``absolute Semple tower'' $(X^a_k,V^a_k)$ associated with
$(X^a_0,V^a_0)=(X,T_X)$. When $X$ is nonsingular, all stages
$(X^a_k,V^a_k)$ of the Semple bundles are nonsingular. This allows to extend
the definition of $(X_k,V_k)$ in case $V$ is singular (but $X$ non singular).
Indeed, let $X'\subset X$ be the Zariski open set on which $V'=V_{|X'}$
is a subbundle of~$T_X$. By functoriality, we get a Semple tower
$(X'_k,V'_k)$ and injections $(X'_k,V'_k)\hookrightarrow(X^a_k,V^a_k)$.
We then simply define $(X_k,V_k)$ to be the pair obtained by taking
the closure of $X'_k$ and $V'_k$ in $X^a_k$ and $V^a_k$ respectively.

\claim 6.8. Proposition {\rm(Direct image formula for invariant
differential operators, [Dem97])}|
Let $E_{k,m}V^*$ be the sheaf of algebraic differential operators 
$f\mapsto P(f_{[k]})$ acting on germs of curves $f:(\bC,T_\bC)\to (X,V)$
such that $P((f\circ\varphi)_{[k]})=
\varphi^{\prime m}P(f_{[k]})\circ\varphi$.
Then
$$
(\pi_{k,0})_*\cO_{X_k}(m)\simeq E_{k,m}V^*.
$$
\endclaim

\claim 6.9. Some tautological morphisms|Let $(X,V)$ be a directed structure
on a nonsingular projective variety~$X$. For every $p=1,\ldots,r$, there
is a tautological morphism
$$
\Phi_{k,p}^{X,V}:\pi_{k,0}^*\Lambda^pV^*\otimes\cO_{X_k}(-(p-1){\ul 1})\to
\Lambda^pV_k^*.
$$
\endclaim

{\it Proof.} First assume that $(X,V)$ is nonsingular. For every integer
$p\geq 1$, the exact sequence (6.3~b) gives
rise to a surjective contraction morphism by the Euler vector field
$$
\Lambda^p \big(\pi_k^*V_{k-1}\otimes\cO_{X_k}(1)\big)^*\to
\Lambda^{p-1}T^*_{X_k/X_{k-1}},
\leqno(6.10)
$$
while (6.3~a) and the dual sequence
$0\to \cO_{X_k}(1)\to V_k^*\to T^*_{X_k/X_{k-1}}\to 0$ induce
an injective wedge multiplication morphism
$$
\Lambda^{p-1}T^*_{X_k/X_{k-1}}\otimes\cO_{X_k}(1)\to
\Lambda^pV_k^*.
\leqno(6.11)
$$
A composition of these two arrows yields a tautological morphism
$$
\pi_k^*\Lambda^pV_{k-1}^*\otimes\cO_{X_k}(-(p-1))\to \Lambda^pV_k^*,
\leqno(6.12)
$$
which turns out to be dual to the isomorphism (6.6) in case  $p=r=\rank V$.
There also exists a (more naive) composed morphism
$$
\pi_k^*\Lambda^pV_{k-1}\otimes\cO_{X_k}(p)\to
\Lambda^pT_{X_k/X_{k-1}}\to \Lambda^pV_k
$$
of dual
$$
\Lambda^pV_k^*\to \pi_k^*\Lambda^pV^*_{k-1}\otimes\cO_{X_k}(-p),
\leqno(6.12')
$$
but it is less interesting for us (the composition of (6.12) and $(6.12')$
is equal to $0$, and $(6.12')$ vanishes for $p=r\,$; these morphisms can be
seen to be dual modulo exchange of $p$ and $r-p$ and contraction with
the determinants). We can iterate (6.12) for all stages
of the tower, and get in this way the desired tautological morphism
$$
\Phi_{k,p}^{X,V}:
\pi_{k,0}^*\Lambda^pV^*\otimes\cO_{X_k}(-(p-1)\ul 1)\to \Lambda^pV_k^*,
$$
By what we have seen, there is in particular an absolute tautological
morphism $\Phi_{k,p}^{X^a,V^a}$ defined at the level of the absolute
Semple tower $(X^a_k,V^a_k)$, and it is clear that $\Phi_{k,p}^{X,V}$
is obtained by restricting $\Phi_{k,p}^{X^a,V^a}$ to $(X_k,V_k)$.
Therefore, our morphisms are also well defined in case $V$ is
singular, by taking the Zariski closure of the regular part in
$(X^a_k,V^a_k)$.\qed
\medskip

Now, let $Z$ be an irreducible algebraic subset of some Semple $k$-jet
bundle $X_k$ over~$X$ ($k$~being arbitrary). We define an {\it induced
directed structure $(Z,W)\hookrightarrow(X_k,V_k)$} by taking the
linear subspace $W\subset T_Z\subset T_{X_k|Z}$ to be the closure of
$T_{Z'}\cap V_k$ taken on a suitable Zariski open set $Z'\subset
Z_{\rm reg}$ where the intersection has constant rank and is a
subbundle of~$T_{Z'}$. Alternatively, one could also take $W$ to be
the closure of $T_{Z'}\cap V_k$ in the absolute $k$-th bundle $V^a_k$.
This produces an {\it induced directed subvariety}
$$
(Z,W)\subset(X_k,V_k)\subset (X^a_k,V^a_k).
\leqno(6.13)
$$
{\bf 6.14. Observation.} One can see that the hypotheses
$Z\subsetneq X_k$ and $\pi_{k,k-1}(Z)=X_{k-1}$
imply $\rank W<\rank V_k$. Otherwise, we would have
$W=(V_k)_{|Z'}\supset T_{X_k/X_{k-1}|Z'}$ over a Zariski open, hence
$Z$ would contain a dense open subset of $X_k$, contradiction.
\medskip

{\bf 6.15. Induced tautological morphisms.} Let $Z$ be an irreducible
subvariety of $X_k$ such that $\pi_{k,0}(Z)=X$, and $(Z,W_k)$ the induced
directed structure. Consider the vertical part $W^v_k\subset W_k$
which is defined to be $W^v_k=\overline{(W_k\cap T_{X_k/X_{k-1}})_{|Z'}}$
where $Z'$ is the Zariski open set where the intersection has
minimal rank (since $\pi_{k,0}(Z)=X$,
we can assume $Z'\subset\pi_{k,0}^{-1}(X')$, where $X'$ is a Zariski open set
where $V_{|X'}$ is a subbundle of $T_X$). Then $W^v_k$ has corank at most $1$
in $W_k$ (since $T_{X_k/X_{k-1}}$ has corank~$1$ in $V_k$). We define
$W_{k-1}\subset \pi_{k,k-1}^*V_{k-1|Z}$ to be the closure of the
preimage of $(W_k^v\otimes \cO_{X_k}(-1))_{|Z'}$ by the restriction to $Z'$
of the morphism
$$
\pi_{k,k-1}^*V_{k-1}\to T_{X_k/X_{k-1}}\otimes \cO_{X_k}(-1)
$$
induced by (6.3~b). Notice that $\rank W_{k-1}=\rank W_k^v+1$ is equal to
$\rank W_k$ or \hbox{$\rank W_k+1$}. We then take
$$
W_{k-1}^v=\overline{(W_{k-1}\cap\pi_{k,k-1}^*T_{X_{k-1}/X_{k-2}})_{|Z''}}
$$
for a suitable Zariski open set $Z''\subset Z$, and obtain in this way
a preimage $W_{k-2}\subset \pi_{k,k-2}^*V_{k-2}$ of
$W_{k-1}^v\otimes\pi_{k,k-1}^*\cO_{X_{k-1}}(-1)$. Inductively, we get a linear
subspace $W_\ell\subset\pi_{k,\ell}^*V_\ell$, and finally
$W_0\subset\pi_{k,0}^*V_0=\pi_{k,0}^*V$,
with
$$
\rank W_0\geq \rank W_1\geq \ldots\geq \rank W_{k-1}\geq \rank W_k.
\leqno(6.15~{\rm a})
$$
Set $r_\ell=\rank W_\ell$. The contraction by the Euler vector field
(6.10) induces by restriction a generically surjective morphism
$$
\Lambda^{p}W_{k-1}^*\otimes\cO_{X_k}(-p)_{|Z}\to
\Lambda^{p-1}(W^v_k)^*.
$$
Now, we also have a sequence $\cO_{X_k}(1)\to W_k^*\to (W_k^v)^*$, which
either gives an isomorphism $\Lambda^{p-1}W_k^*=
\Lambda^{p-1}(W^v_k)^*$ when $W_k^v=W_k$ (in that case, the arrow
$\cO_{X_k}(1)\to W_k^*$ is equal to~$0$), or a generically surjective
morphism
$$
\Lambda^{p-1}(W^v_k)^*\otimes\cO_{X_k}(1)\to\Lambda^pW_k^*
$$
induced by (6.11) when $\rank W_k^v=\rank W_k-1$. Therefore we obtain
by composition generically surjective morphisms
$$
\eqalign{
&\Lambda^{p}W_{k-1}^*\otimes\cO_{X_k}(-p)_{|Z}\to
\Lambda^{p-1}W_k^*\quad\hbox{when $r_{k-1}=r_k+1$},\cr
&\Lambda^{p}W_{k-1}^*\otimes\cO_{X_k}(1-p)_{|Z}\to
\Lambda^pW_k^*\quad\hbox{when $r_{k-1}=r_k$}.\cr}
$$
Replacing $p$ by $r_{k-1}-p$, we obtain in both cases a generically surjective
morphism
$$
\Lambda^{r_{k-1}-p}W_{k-1}^*\to
\Lambda^{r_k-p}W_k^*\otimes\cO_{X_k}(2r_{k-1}-r_k-p-1)_{|Z},
$$
and inductively, we obtain for every $p\geq 0$ a generically
surjective morphism
$$
\Lambda^{r_0-p}W_0^*\to \Lambda^{r_k-p}W_k^*\otimes
\cO_{X_k}(\ul a-p\,\ul 1)_{|Z},\quad
\ul a=(a_\ell),~~a_\ell=2r_{\ell-1}-r_\ell-1.
\leqno(6.15~{\rm b})
$$
For $p=0$, we get in particular a generically surjective morphism
$$
\Lambda^{r_0}W_0^*\to \Lambda^{r_k}W_k^*\otimes
\cO_{X_k}(\ul a)_{|Z},\quad
\leqno(6.15~{\rm c})
$$
which is especially interesting since we are in rank~$1$. Since
$W_0\subset(\pi_{k,0}^* V)_{|Z}\subset(\pi_{k,0}^* T_X)_{|Z}$, we also have
generically surjective morphisms
$$
(\pi_{k,0}^*\Lambda^{r_0}T_X^*)_{|Z}\to
(\pi_{k,0}^*\Lambda^{r_0}V^*)_{|Z}\to\Lambda^{r_0}W_0^*.
$$
Also, all of these morphisms are obtained by taking either restrictions
of forms or contractions by the Euler vector fields of the various stages.
With respect to smooth hermitian metrics on the (nonsingular) absolute
Semple tower $(X^a_k,V^a_k)$, bounded forms are certainly mapped to
bounded forms. From this we conclude~:

\claim 6.16. Corollary|Let $X$ be a nonsingular variety, $(X,V)$ a directed
structure, and \hbox{$Z\subset X_k$} an irreducible subvariety
in the $k$-th stage of the Semple tower $(X_k,V_k)$ such that
 \hbox{$\pi_{k,0}(Z)=X$}. Then the induced directed variety $(Z,W)$ has the
following property: there exists $r_0\geq r_k:=\rank W$, a weight
$\ul a\in\bN^k$ and a generically surjective sheaf morphism
$$
(\pi_{k,0}^*\Lambda^{r_0}T_X^*)_{|Z}\to(\pi_{k,0}^*{}^b\Lambda^{r_0}V^*)_{|Z}\to
{}^b\cK_W\otimes\cO_{X_k}(\ul a)_{|Z}.
$$
\endclaim
\medskip

{\bigbf 7. Tautological morphisms and differential equations}

The tautological morphisms give a potential technique for controlling
inductively the positivity of the canonical sheaf ${}^b\cK_W$ for all induced
directed structures $(Z,W)$. We rely on the following
statement that was already observed in [Dem14].

\claim 7.1. Proposition|Let $X$ be a nonsingular variety, $(X,V)$ a directed
structure, and $(Z,W)$ the induced directed structure on an irreducible
subvariety $Z\subset X_k$. Assume that there exists a weight
$\ul a\in\bQ_+^k$ such that ${}^b\cK_W\otimes\cO_{X_k}(\ul a)_{|Z}$ is big.
Then there exists \hbox{$Y\subsetneq Z$} and finitely many irreducible
directed subvarieties $(Z'_{\ell,\alpha},W'_{\ell,\alpha})\subset
(X_{k+\ell},V_{k+\ell})$ contained in the Semple tower
$(Z_\ell,W_\ell)$ of $(Z,W)$, with \hbox{$\rank W'_{\ell,\alpha}<\rank W$},
such that
every entire curve $f:(\bC,T_\bC)\to(Z,W)$ is either contained in $Y$, or
has a lifting $f_{[\ell]}:(\bC,T_\bC)\to (Z'_{\ell,\alpha},W'_{\ell,\alpha})$
for some $\alpha$.
\endclaim
\medskip

{\it Idea of the proof.} One can take $Y\subset Z$ to be a subvariety such that
$Z\ssm Y\subset Z_{\rm reg}$, $W_{|Z\ssm Y}$ is non singular, and for some
ample divisor $A'$ on $Z$ and $m_0\gg 1$, the sheaf
$({}^b\cK_W\otimes\cO_{X_k}(\ul a)_{|Z})^{\otimes m_0}\otimes\cO_{Z}(-A')$ is
generated by sections over $Z\ssm Y$. The assumption
on ${}^b\cK_W\otimes\cO_{X_k}(\ul a)_{|Z}$ allows to use again holomorphic
Morse inequalities to construct a nontrivial section 
$$
\sigma\in H^0\big((\tilde Z_\ell,\cO_{\tilde Z_\ell}(\ul b)\otimes
\tilde\pi_{\ell,0}^*(\mu^*\cO_{X_k}(p\,\ul a)_{|Z}\otimes
\cO_{\tilde Z}(-A))\big)
$$
on some nonsingular modification $\mu:\tilde Z\to Z$.
This method produces a new differential equation $\sigma=0$ that yields the
desired subvariety $Z'_\ell=\bigcup Z'_{\ell,\alpha}\subsetneq Z_\ell$.\qed
\medskip

We will also make use of the following observation due
to Laytimi and Nahm [LNa99], valid for any
holomorphic vector bundle $E\to X$ on a projective manifold
(one can more generally consider torsion free coherent sheaves).

\claim 7.2. Definition|We say that $E\to X$ {\rm is strongly big} if
for any $A\in\Pic(X)$ ample, the symmetric powers
$S^mE\otimes\cO(-A)$ are generated by their sections over
a Zariski open set of $X$, for a sufficiently large integer $m\gg 1$.
\endclaim

\claim 7.3. Lemma|For every $p\ge 0$, if $\Lambda^pE$ is big, then
$\Lambda^{p+1}E$ is also big.
\endclaim
\medskip

{\it Sketch of proof.} This is a consequence of the fact,
observed by [LNa99, Lemma 2.4],
that for $k\gg 1$, the Schur irreducible components of
$(\Lambda^{k+1}E)^{\otimes kp}$ all appear in $(\Lambda^kE)^{(k+1)p}$.
We thank L.~Manivel for pointing out this simple argument.
\qed

\medskip

From this, one can now derive the following statement. 

\claim 7.4. Theorem|
Let $(X,V)$ be a directed variety. Assume that ${}^b\Lambda^pV^*$
is a strongly big sheaf for some $p\in\bN^*$, $p\leq r=\rank V$.
{\parindent=6.5mm
\item{\rm(a)} If $p=1$, $(X,V)$ satisfies the generalized GGL conjecture, i.e.,
there exists a subvariety $Y\subsetneq X$ containing all entire curves
$f:(\bC,T_\bC)\to(X,V)$.
\item{\rm(b)} If $p\ge 2$, there exists a subvariety $Y\subsetneq X$, an
integer $k\in \bN$ and a
finite collection of induced directed subvarieties $(Z_\alpha,W_\alpha)
\subset (X_k,V_k)$ with $Z_\alpha$ irreducible, \hbox{$\rank W_\alpha\leq p-1$},
such that all entire curves $f:(\bC,T_\bC)\to(X,V)$ satisfy either
$f(\bC)\subset Y$ or have a $k$-jet lifting
$f_{[k]}:(\bC,T_\bC)\to(Z_\alpha,W_\alpha)$ for some
$\alpha$.
\item{\rm(c)} In particular, if $p=2$, all entire curves $f:(\bC,T_\bC)\to(X,V)$
are either contained in $Y\subsetneq X$, or they are
tangent to a {\rm rank~$1$ foliation} on a subvariety $Z=\bigcup Z_\alpha
\subset X_k$. This
implies that the latter~curves are supported by the parabolic leaves of
the above foliations, which can be parametrized as a subspace of a
finite dimensional variety.
\item{\rm(d)} The subvariety $Y$ described in {\rm(a)}, {\rm(b)} or {\rm(c)}
can be taken to be any subvariety such that
$S^{m_0}({}^b\Lambda^pV^*)\otimes\cO_X(-A)$ is generated by sections over
$X\ssm Y$, for a suitable $A\in\Pic(X)$ ample and $m_0\gg 1$. In particular, if
${}^b\Lambda^pV^*$ is ample, one can take $Y=\emptyset$.
\vskip0pt}
\endclaim
\medskip

{\it Proof.} (a) As we already observed, the rank 1 case is an easy
consequence of the Ahlfors-Schwarz lemma. Also, (c) is a particular case
of~(b), so we only have to check (b) and~(d).

(b) For $p=r$, the statement is a consequence of Corollary 4.12.
In general, we decompose all occurring subvarieties into irreductible 
components $(Z_\alpha,W_\alpha)$ and apply descending induction on
$r_\alpha=\rank W_\alpha$ for each of them. As long as $r_\alpha\geq p$, 
The assumption on $V$ combined with Corollary 6.16 implies that
there exists a weight $\ul a\in\bQ_+^k$ such that
${}^b\cK_{W_\alpha}\otimes\cO_{X_k}(\ul a)_{|Z_\alpha}$ is big (for this we use
the fact that $\cO_{X_k}(1)$ is relatively big with respect
to $X_k\to X$). Then Proposition 7.1 either produces a
subvariety $Y'\subsetneq Z_\alpha$ (in which case we consider the irreducible
components $Y'_\beta$ and apply descending induction on
$\dim Y'_\beta$, if the rank of the induced structure does not
decrease right away), or we get irreducible directed subvarieties
$(Z'_{\ell,\beta},W'_{\ell,\beta})$ in the $\ell$-th stage of the Semple tower
of $(Z_\alpha,W_\alpha)$, with $r'_{\ell,\beta}=\rank W'_{\ell,\beta}<r_\alpha$.
The induction hypothesis applies as long as
$p\leq r'_{\ell,\beta}<r_\alpha$. Therefore we end up with
$r'_{\ell,\beta}<p$ after finitely many iterations.

(d) is a consequence of the technique of proof of [Dem11] and [Dem14]. In fact
we use singular hermitian metrics on $V$ satisfying suitable positivity
properties, and $Y$ can be taken to be their set of poles.\qed
\medskip

{\bf 7.5. Remark.} It would be interesting to know if the rather restrictive
bigness hypothesis for ${}^b\Lambda^pV^*$ can be relaxed, assuming instead
suitable semistability conditions. In fact, if none of the algebraic
hypersurfaces contains all of the entire curves, one can show that there exists
an Ahlfors current that defines a mobile bidegree $(n-1,n-1)$ class, and
one could try to use the Harder-Narasimhan filtration of $V$ with respect to
such mobile classes.
\medskip

{\bf 7.6. Logarithmic and orbifold directed versions.}
More generally, let $\Delta=\sum\Delta_j$ be a reduced
normal crossing divisor in $X$. We want to study entire
curves $f:\bC\to X\ssm\Delta$ drawn in the complement of $\Delta$.
At a point where $\Delta=\{z_1\ldots z_p=0\}$ one defines the
{\it logarithmic cotangent sheaf $T^*_{X\langle\Delta\rangle}$} to
be generated by ${dz_1\over z_1},...,{dz_p\over z_p},dz_{p+1},...,dz_n$.
The results stated above can easily be extended to the logarithmic case.
In particular, we obtain the following statement.

\claim 7.7. Theorem|
If $\Lambda^2T^*_{X\langle\Delta\rangle}$ is strongly big on $X$, there
exists a subvariety $Y\subsetneq X$ and rank 1 foliations $\cF_\alpha$
on some subvarieties $Z_\alpha\subset X_k$ of the $k$-jet bundle,
such that all entire curves $f:\bC\to X\ssm\Delta$ are either contained
in $Y$ or have a $k$-jet lifting $f_{[k]}$ that is contained in
some $Z_\alpha$ and tangent to $\cF_\alpha$. When
$\Lambda^2T^*_{X\langle\Delta\rangle}$ is ample, we can take $Y=\emptyset$.
\endclaim
\medskip

One can obtain even more general versions dealing with entire
curves $f:(\bC,T_\bC)\to (X,V)$ that are tangent to $V$ and avoid a
normal crossing divisor $\Delta$ transverse to $V$
({\it logarithmic~case}), or meet $\Delta=
\sum(1-{1\over\rho_j})\Delta_j$ with multiplicities${}\geq \rho_j$
along $\Delta_j$ ({\it orbifold case}).
Such statements are the subject of a work in progress with
F.\ Campana, L.\ Darondeau and E.\ Rousseau. In this setting, the
positivity conditions have to be expressed not just for
the orbifold cotangent bundle, but also for the ``{\it derived
orbifold cotangent bundles}'' of higher order
$V^*\langle\Delta^{(s)}\rangle$, associated with the divisors
$\Delta^{(s)}=\sum_j\big(1-{s\over\rho_j}\big)_+\Delta_j$, $1\leq s\leq k$.
\bigskip\bigskip

\centerline{\bigbf References}
\medskip

\bibitem[Bon93]&Bonavero, L.:&In\'egalit\'es de Morse holomorphes
singuli\`eres.&C.~R.\ Acad.\ Sci.\ Paris S\'er.~I Math.\ {\bf 317} (1993)
1163-–1166&

\bibitem[Dem85]&Demailly, J.-P.:&Champs magn\'etiques et in\'egalit\'es de
Morse pour la $d''$-coho\-mo\-logie.&Ann.\ Inst.\ Fourier
(Grenoble) {\bf 35} (1985) 189--229&

\bibitem[Dem97]&Demailly, J.-P.:&Algebraic criteria for Kobayashi
hyperbolic projective varieties and jet differentials.&AMS Summer
School on Algebraic Geometry, Santa Cruz 1995, Proc.\ Symposia in
Pure Math., vol.\ {\bf 062.2}, ed.\ by J.~Koll\'ar and R.~Lazarsfeld,
Amer.\ Math.\ Soc., Providence, RI (1997), 285–-360&

\bibitem[Dem11]&Demailly, J.-P.:&Holomorphic
Morse Inequalities and the Green-Griffiths-Lang Conjecture.&Pure
and Applied Math.\ Quarterly {\bf 7} (2011), 1165--1208&

\bibitem[Dem14]&Demailly, J.-P.:&Towards
the Green-Griffiths-Lang conjecture.&Conference\break
``Analysis and Geometry'', Tunis,
March 2014, in honor of Mohammed Salah Baouendi, 
ed.\ by A.~Baklouti, A.\ El Kacimi, S.~Kallel, N.~Mir, Springer
Proc.\ Math.\ Stat., {\bf 127}, Springer-Verlag, 2015,
141--159&

\bibitem[Dem20]&Demailly, J.-P.:&Recent results on the Kobayashi
and Green-Griffiths-Lang conjectures&Japanese J.~of Math., {\bf 15}
(2020), 1--120&

\bibitem[Ete19]&Etesse, A.:&Ampleness of Schur powers of cotangent bundles
and $k$-hyperbolicity&Res.\ Math.\ Sci., vol.~8, article \#~7 (2021)&

\bibitem[GGr80]&Green, M., Griffiths, P.A.:&Two applications of algebraic
geometry to entire holomorphic mappings.&The Chern Symposium 1979,
Proc.\ Internal.\ Sympos.\ Berkeley, CA, 1979, Springer-Verlag, New York
(1980), 41--74&

\bibitem[Lan86]&Lang, S.:& Hyperbolic and Diophantine analysis.&Bull.\
Amer.\ Math.\ Soc.\ {\bf 14} (1986), 159--205&

\bibitem[Lan87]&Lang, S.:& Introduction to complex hyperbolic
spaces.&Springer-Verlag, New York (1987)&

\bibitem[LNa99]&Laytimi, F., Nahm, W.:&Vanishing theorems for products of
exterior and symmetric powers&arXiv:math/9809064v2 [math.AG], 24 Feb 1999&

\bibitem[McQ98]&McQuillan M.:&Diophantine approximations and foliations&Publ.\
Math.\ I.H.\'E.S.\ {\bf 87} (1998) 121--174&

\bibitem [Sem54]&Semple, J.G.:&Some investigations in the geometry of
curves and surface elements.&Proc.\ London Math.\ Soc.\ (3) {\bf 4}
(1954), 24--49&

\bibitem[SYe97]&Siu, Y.T., Yeung, S.K.:&Defects for ample divisors of
Abelian varieties, Schwarz lemma and hyperbolic hypersurfaces of low
degree&Amer.~J.\ Math.\ {\bf 119} (1997), 1139--1172&
\bigskip
Jean-Pierre Demailly\\
Institut Fourier, Universit\'e Grenoble Alpes\\
100, rue des Maths, 38610 Gi\`eres, France\\
{\it e-mail}\/: jean-pierre.demailly@univ-grenoble-alpes.fr
\bigskip
(June 11, 2021; printed on \today)
\bye